\documentclass[11pt,reqno]{amsart}
\usepackage{amssymb,amscd,amsbsy}
\usepackage{amsthm}
\newcommand{\lb}{\linebreak}

\renewcommand{\a}{\alpha}

\renewcommand{\d}{\delta}
\newcommand{\e}{\varepsilon}

\newcommand{\z}{\zeta}

\renewcommand{\l}{\lambda}

\newcommand{\s}{\sigma}
\renewcommand{\t}{\tau}

\newcommand{\f}{\varphi}
\renewcommand{\o}{\omega}

\newcommand{\D}{\Delta}
\renewcommand{\L}{\Lambda}

\renewcommand{\O}{\Omega}

\newcommand{\n}{\nabla}

\newcommand{\B}{{\mathcal B}}

\newcommand{\E}{{\mathcal E}}
\newcommand{\cd}{{\mathcal D}}
\newcommand{\F}{{\mathcal F}}

\newcommand{\h}{{\mathcal H}}

\newcommand{\K}{{\mathcal K}}
\newcommand{\cL}{{\mathcal L}}

\newcommand{\X}{{\mathcal X}}
\newcommand{\Y}{{\mathcal Y}}

\newcommand{\T}{{\Bbb T}}
\newcommand{\pp}{{\Bbb P}}
\newcommand{\dd}{{\Bbb D}}
\newcommand{\R}{{\Bbb R}}
\newcommand{\Z}{{\Bbb Z}}

\newcommand{\0}{{\boldsymbol{0}}}

\newcommand{\bs}{\boldsymbol}

\newcommand{\m}{{\boldsymbol m}}
\newcommand{\bS}{{\boldsymbol S}}

\newcommand{\rf}[1]{(\ref{#1})}

\newcommand{\df}{\stackrel{\mathrm{def}}{=}}

\newcommand{\spn}{\operatorname{span}}
\newcommand{\supp}{\operatorname{supp}}
\newcommand{\clos}{\operatorname{clos}}

\newcommand{\const}{\operatorname{const}}

\newcommand{\eeq}{\end{equation}}
\newcommand{\beq}{\begin{equation}}
\newcommand{\bay}{\begin{eqnarray}}
\newcommand{\ba}{\begin{align*}}
\newcommand{\ea}{\end{align*}}
\newcommand{\ey}{\end{eqnarray}}
\newcommand{\bey}{\begin{eqnarray*}}
\newcommand{\eey}{\end{eqnarray*}}

\newcommand{\imp}{\Rightarrow}
\newcommand{\be}{\infty}

\newcommand{\bl}{\blacksquare}

\newcommand{\Pf}{{\bf Proof. }}

\newcommand{\ov}{\overline}

\newtheorem{thm}{\hspace{\parindent}Theorem}[section]

\newtheorem{cor}[thm]{\hspace{\parindent}Corollary}
\newtheorem{lem}[thm]{\hspace{\parindent}Lemma}

\theoremstyle{remark}

\newtheorem*{rem*}{Remark}

\newcommand\up{\upsilon}

\newcommand\fM{\frak M}
\newcommand\cZ{\mathcal{Z}}
\newcommand\dg{\frak D}

\begin{document}

\newcommand{\vse}{\vspace{.2in}}
\numberwithin{equation}{section}

\title{ Differentiability of functions of contractions}
\author{V.V. Peller}
\thanks{The author is partially supported by NSF grant DMS 0700995}
\maketitle

\newcommand{\mt}{{\mathcal T}}

\begin{abstract}
In this paper we study differentiability properties of the map $T\mapsto\f(T)$, where $\f$ is a given function in the disk-algebra
and $T$ ranges over the set of contractions on Hilbert space. We obtain sharp conditions (in terms of Besov spaces) for differentiability and existence of higher derivatives. We also find explicit formulae for directional derivatives (and higher derivatives) in terms of double (and multiple) operator integrals with respect to semi-spectral measures.
\end{abstract}

\setcounter{section}{0}
\section{\bf Introduction}
\setcounter{equation}{0}

\

The purpose of this paper is to study differentiability properties of functions
$$
T\mapsto\f(T),
$$
for a given function $\f$ analytic in the unit open disk $\dd$ and continuous in the closed disk (in other words $\f$ belongs to the disk-algebra $C_A$), where $T$ ranges over the set of contractions (i.e., operators of norm at most $1$) on Hilbert space.

Recall that by von Neumann's inequality,
\bay
\label{vN}
\|\f(T)\|\le\max_{|\z|\le1}|\f(\z)|
\ey
for an arbitrary contraction $T$ on Hilbert space and an arbitrary polynomial $\f$. This allows one to define a functional calculus
$$
\f\mapsto\f(T), \quad\f\in C_A,
$$
for an arbitrary contraction $T$. Moreover, for this functional calculus von Neumann's inequality holds \rf{vN} holds.

For contractions $T$ and $R$ on Hilbert space, we consider the one-parameter family of contractions
$$
T_t=(1-t)T+tR,\quad\0\le t\le1,
$$
and we study differentiability properties of the map
\bay
\label{opf}
t\mapsto\f(T_t)
\ey
for a given function $\f$ in $C_A$.

The study of the problem of differentiability of functions of self-adjoint operators on Hilbert space was initiated 
By Daletskii and S.G. Krein in \cite{DK}. They showed that for a function $f$ on the real line $\R$ of class $C^2$ and for 
bounded self-adjoint operators $A$ and $B$ the function
\bay
\label{sso}
t\mapsto f(A+tB)
\ey
is differentiable in the operator norm and the derivative can be computed in terms of double operator integrals:
\bay
\label{DKF}
\frac{d}{dt}f(A+tB)\Big|_{t=0}=\iint\limits_{\R\times\R}\frac{f(x)-f(y)}{x-y}\,dE_A(x)\,B\,dE_A(y),
\ey
where $E_A$ is the spectral measure of $A$. The expression on the right is a double operator integral. The beautiful theory of double operator integrals due to Birman and Solomyak was created later in \cite{BS1}, \cite{BS2}, and \cite{BS3} (see also 
the survey article \cite{BS4}).

The condition $f\in C^2$ was relaxed by Birman and Solomyak in  \cite{BS3}: they proved that the function \rf{sso} is differentiable and the Daletskii--Krein formula \rf{DKF} holds under the condition that $f$ is differentiable and the derivative $f'$ satisfies a H\"older condition of order $\a$ for some $\a>0$. The approach of Birman and Solomyak is based on their formula
\bay
\label{BSdof}
f(A+B)-f(A)=\iint\limits_{\R\times\R}\frac{f(x)-f(y)}{x-y}\,dE_{A+B}(x)\,B\,dE_A(y).
\ey
Actually, Birman and Solomyak showed in \cite{BS3} that formula \rf{BSdof}
is valid under the condition that the divided difference $\cd f$ is a Schur multiplier of the space of all bounded linear operators
(see \S\,2.2 for the definitions).

However, it follows from the results of Farforovskaya in \cite{F} that the condition $f\in C^1$ is not sufficient for the differentiability of the map \rf{sso}.

A further improvement was obtained in \cite{Pe1}: it was shown that the the function \rf{sso} is differentiable and
\rf{DKF} holds under the assumption that $f$ belongs to the Besov space $B_{\be1}^1(\R)$ (see subsection 2.5). Moreover, in \cite{Pe1} a necessary condition was also found: $f$ must locally belong to the Besov space $B_1^1(\R)=B^1_{11}(\R)$. This necessary condition also implies that the condition $f\in C^1$ is not sufficient. Actually, in \cite{Pe1} a stronger necessary condition was also obtained, see \S\,2.3 for a further discussion. Finally, we mention another sufficient condition obtained in \cite{ABF} which is slightly better than the condition $f\in B_{\be1}^1(\R)$, though I believe it is more convenient to work with Besov spaces. 

Note that similar results were obtained in \cite{Pe4} in the case when $A$ is an unbounded self-adjoint operator and $B$ is a bounded self-adjoint operator.

The problem of the existence of higher derivatives of the function \rf{sso} was studied in \cite{S} where it was shown that under certain assumptions on $f$, the function \rf{sso} has second derivative that can be expressed in terms of 
the following triple operator integral:
$$
\frac{d^2}{dt^2}f(A+tB)\Big|_{t=0}=\iiint\limits_{\R\times\R\times\R}\left(\cd^2\f\right)(x,y,z)\,dE_A(x)\,B\,dE_A(y)\,B\,dE_A(z),
$$
where $\cd^2\f$ stands for the divided difference of order 2 (see \S\,2 for the definition).
To interpret triple operator integrals, repeated integration was used in \cite{S}. However, the class of integrable functions in \cite{S} was rather narrow and the assumption on $f$ imposed in \cite{S} for the existence of the second operator derivative was too restrictive.
Similar results are also obtained in \cite{S} for the $n$th derivative and multiple operator integrals.

In \cite{Pe8} a new approach to multiple operator integrals was given. It is based on integral projective tensor products of $L^\be$ spaces and gives a much broader class of integrable functions than under the approach of \cite{S}. It was shown in \cite{Pe8} that under the assumption that $f$ belongs to the Besov space $B_{\be1}^n(\R)$ the function \rf{sso} has $n$ derivatives and the $n$th derivative can be expressed in terms of a multiple operator integral. Similar results were also obtained in \cite{Pe8} in the case of an unbounded self-adjoint operator $A$.

Note that Besov spaces $B_{\be1}^n(\R)$ arise in operator theory on many different occasions, see \cite{Pe},
 \cite{Pe3}, \cite{Pe5},
\cite{Pe7}.

Let us mention here another formula by Birman and Solomyak (see \cite{BS4}) for commutators. Suppose that $A$ is a self-adjoint operator and $Q$ is a bounded linear operator. Then 
\bay
\label{BScom}
f(A)Q-Qf(A)=\iint\limits_{\R\times\R}\frac{f(x)-f(y)}{x-y}\,dE_A(x)\,(AB-BA)\,dE_A(y),
\ey
which is also valid under the assumption that $\cd f$ is a Schur multiplier of the space of all bounded linear operators.

To study the problem of differentiability of functions of unitary operators, we should consider a Borel function $f$ on the unit circle
$\T$ and the map
$$
U\mapsto f(U),
$$
where $U$ is a unitary operator on Hilbert space. If $U$ and $V$ are unitary operators and $V=e^{{\rm i}A}U$, 
where $A$ is a self-adjoint operator, we can consider the one-parametric family of unitary operators
$$
e^{{\rm i}tA}U,\quad0\le t\le1,
$$
and study the question of the differentiability of the function 
$$
t\mapsto \big(e^{{\rm i}tA}U\big)
$$
and the question of the existence of its higher derivatives.
The results in the case of unitary operators are similar to the results for self-adjoint operators, see \cite{BS3}, \cite{Pe1},
\cite{ABF}, \cite{Pe8}.

In this paper we study the case of functions of contractions. This study was initiated in \cite{Pe2}, where the Lipschitz property was
studied. Recently, in \cite{KS} new results on operator Lipschitz functions of contractions were obtained, see \S~2.3 for more detailed information.

It turns out that the right tool to study differentiability properties of functions of contractions is double (and multiple) operator integrals with respect to semi-spectral measures. Note that even if both contractions $T$ and $R$ are unitary operators, the differentiability problem for this pair is different from the differentiability problem for unitary operators.

In \S\,3 we define double and multiple operator integrals with respect to semi-spectral measures.
In \S\,4 we obtain an analog of the Birman--Solomyak formulae \rf{BSdof} and \rf{BScom} for semi-spectral measures. 
Then we obtain in \S\,5 conditions on a function $\f\in C_A$ for the differentiability of the map \rf{opf} in the operator norm
as well as conditions for the existence of higher operator derivatives. We also obtain in \S\,5 formuale for the derivatives of 
\rf{opf} in terms of multiple operator integrals with respect to semi-spectral measures.
Finally, in \S\,6 we study the problem of differentiability of the function \rf{opf} in the Hilbert--Schmidt norm. 

In \S\,2 we give necessary information on Besov spaces, double and multiple operator integrals, and semi-spectral measures.

I would like to express my gratitude to Victor Shulman for stimulating discussions.

\

\section{\bf Preliminaries}
\setcounter{equation}{0}

\

We are going to collect in this section necessary information on Besov spaces, double operator integrals, multiple operator itegrals, and semi-spectral measures.

{\bf2.1. Besov spaces.}  Let $1\le p,\,q\le\be$ and $s\in\R$. The Besov class $B^s_{pq}$ of functions (or
distributions) on $\T$ can be defined in the following way. Let $w$ be a piecewise linear function on $\R$ such
that
$$
w\ge0,\quad\supp w\subset\left[\frac12,2\right],\quad w(1)=1,
$$
and $w$ is a linear function on the intervals $[1/2,1]$ and $[1,2]$.

Consider the trigonometric polynomials $W_n$, and $W_n^\#$ defined by
$$
W_n(z)=\sum_{k\in\Z}w\left(\frac{k}{2^n}\right)z^k,\quad n\ge1,\quad W_0(z)=\bar z+1+z,\quad
\mbox{and}\quad W_n^\#(z)=\ov{W_n(z)},\quad n\ge0.
$$
Then for each distribution $\f$ on $\T$,
$$
\f=\sum_{n\ge0}\f*W_n+\sum_{n\ge1}\f*W^\#_n.
$$
The Besov class $B^s_{pq}$ consists of functions (in the case $s>0$) or distributions $\f$ on $\T$
such that
$$
\big\{\|2^{ns}\f*W_n\|_{L^p}\big\}_{n\ge0}\in\ell^q\quad\mbox{and}
\quad\big\{\|2^{ns}\f*W^\#_n\|_{L^p}\big\}_{n\ge1}\in\ell^q
$$
Besov classes admit many other descriptions. In particular, for $s>0$, the space $B^s_{pq}$ admits the
following characterization. A function $\f$ belongs to $B^s_{pq}$, $s>0$, if and only if
$$
\int_\T\frac{\|\D^n_\t f\|_{L^p}^q}{|1-\t|^{1+sq}}d\m(\t)<\be\quad\mbox{for}\quad q<\be
$$
and
$$
\sup_{\t\ne1}\frac{\|\D^n_\t f\|_{L^p}}{|1-\t|^s}<\be\quad\mbox{for}\quad q=\be,
$$
where $\m$ is normalized Lebesgue measure on $\T$, $n$ is an integer greater than $s$ and $\D_\t$ is
the difference operator: 
$(\D_\t f)(\z)=f(\t\z)-f(\z)$, $\z\in\T$.

We are going to use the notation $B_p^s$ for $B_{pp}^s$.

It is easy to see from the definition of Besov classes that the Riesz projection $\pp_+$,
$$
\pp_+f=\sum_{n\ge0}\hat f(n)z^n,
$$
is bounded on $B^s_{pq}$ and functions in $\big(B^s_{pq}\big)_+\df\pp_+B^s_{pq}$ admit a natural extension to the unit disk $\dd$, they are analytic in $\dd$ and the functions in $\big(B^s_{pq}\big)_+$ admit the following description:
$$
f\in \big(B^s_{pq}\big)_+\Leftrightarrow
\int_0^1(1-r)^{q(n-s)-1}\|f^{(n)}_r\|^q_p\,dr<\be,\quad q<\be,
$$
and
$$
f\in \big(B^s_{p\be}\big)_+\Leftrightarrow
\sup_{0<r<1}(1-r)^{n-s}\|f^{(n)}_r\|_p<\be,
$$
where $f_r(\z)\df f(r\z)$ and $n$ is a nonnegative integer greater than $s$.

In a similar way one can define (homogeneous) Besov space $B_{pq}^s$ of functions (distributions) on $\R$.

We refer the reader to \cite{Pee} and \cite{Pe6} for more detailed information on Besov spaces.

\medskip

{\bf2.2. Double operator integrals.} In this subsection we give a brief introduction in the theory of double operator integrals
developed by Birman and Solomyak in \cite{BS1}, \cite{BS2}, and \cite{BS3}, see also their survey \cite{BS5}.

Let $(\X_1,E_1)$ and $(\X_2,E_2)$ be spaces with spectral measures $E_1$ and $E_2$
on a Hilbert spaces $\h_1$ and $\h_2$. Let us first define double operator integrals
\bay
\label{doi}
\int\limits_{\X_1}\int\limits_{\X_2}\Phi(\l,\mu)\,d E_1(\l)\,Q\,dE_2(\mu),
\ey
for bounded measurable functions $\psi$ and operators $Q:\h_2\to\h_1$
of Hilbert Schmidt class $\bS_2$. Consider the set function $F$ whose values are orthogonal
projections on the Hilbert space $\bS_2(\h_2,\h_1)$ of Hilbert--Schmidt operators from $\h_1$ to $\h_1$, which is defined 
on measurable rectangles by
$$
F(\L\times\D)T=E_1(\L)QE_2(\D),\quad Q\in\bS_2(\h_2,\h_1),
$$ 
$\L$ and $\D$ being measurable subsets of $\X$ and $\Y$. It was shown in \cite{BS4} that $F$ extends to a spectral measure on 
$\X_1\times\X_2$ and if $\psi$ is a bounded measurable function on $\X_1\times\X_2$, we define
$$
\int\limits_{\X_1}\int\limits_{\X_2}\Phi(\l,\mu)\,d E_1(\l)\,Q\,dE_2(\mu)=
\left(\,\,\int\limits_{\X_1\times\X_2}\Phi\,dF\right)Q.
$$
Clearly,
$$
\left\|\,\,\int\limits_{\X_1}\int\limits_{\X_2}\Phi(\l,\mu)\,dE_1(\l)\,Q\,dE_2(\mu)\right\|_{\bS_2}
\le\|\Phi\|_{L^\be}\|T\|_{\bS_2}.
$$

It is easy to see from the definition of double operator integrals in the case $Q\in\bS_2$ 
that if $\{\Phi_n\}_{n\ge1}$ is a sequence of measurable functions such that 
$$
\lim_{n\to\be}\Phi_n(\l,\mu)=\Phi(\l,\mu),\quad \l\in\X_1,~\mu\in\X_2,
$$
and
$$
\sup_n\sup_{\l,\,\mu\in\X_1\times\X_2}|\Phi_n(\l,\mu)|<\be,
$$
then
\bay
\label{poto}
\lim_{n\to\be}\left\|~\int\limits_{\X_1}\int\limits_{\X_2}\Phi_n(\l,\mu)\,d E_1(\l)\,Q\,dE_2(\mu)
-\int\limits_{\X_1}\int\limits_{\X_2}\Phi(\l,\mu)\,d E_1(\l)\,Q\,dE_2(\mu)\right\|_{\bS_2}=0.
\ey

If the transformer
$$
Q\mapsto\int\limits_{\X_1}\int\limits_{\X_2}\Phi(\l,\mu)\,d E_1(\l)\,Q\,dE_2(\mu)
$$
maps the trace class $\bS_1$ into itself, we say that $\Phi$ is a {\it Schur multiplier of $\bS_1$ associated with 
the spectral measure $E_1$ and $E_2$}. In
this case the transformer
\bay
\label{tra}
Q\mapsto\int\limits_{\X_2}\int\limits_{\X_1}\Phi(\l,\mu)\,d E_2(\mu)\,Q\,dE_1(\l),\quad Q\in \bS_2(\h_1,\h_2),
\ey
extends by duality to a bounded linear transformer on the space of bounded linear operators from $\h_1$ to $\h_2$
and we say that the function $\Psi$ on $\X_2\times\X_1$ defined by 
$$
\Psi(\mu,\l)=\Phi(\l,\mu)
$$
is {\it a Schur multiplier of the space of bounded linear operators associated with $E_2$ and $E_1$}.
We denote the space of such Schur multipliers by $\fM(E_2,E_1)$


In \cite{BS3} it was shown that if $A$ is a self-adjoint operator (not necessarily bounded),
$K$ is a bounded self-adjoint operator and if
$\f$ is a continuously differentiable 
function on $\R$ such that the divided difference $\dg\f$ defined by
$$
(\dg\f)(\l,\mu)=\frac{\f(\l)-\f(\mu)}{\l-\mu}
$$
is a Schur multiplier of 
the space of bounded linear operators with respect to the spectral measures of $A+K$ and $A$, then
\bay
\label{BSF}
\f(A+K)-\f(A)=\iint\limits_{\R\times\R}\frac{\f(\l)-\f(\mu)}{\l-\mu}\,dE_{A+K}(\l)K\,dE_A(\mu)
\ey
and
$$
\|\f(A+K)-\f(A)\|\le\const\|\f\|_{\fM}\|K\|,
$$
where $\|\f\|_{\fM}$ is the norm of $\f$ in $\fM(E_{A+K},E_A)$. 

The same formula \rf{BSF} holds in the case $K$ is a Hilbert--Schmidt operator and $\f$ is a Lipschitz function (in this case $\dg\f$ is not necessarily defined on the diagonal of $\R\times\R$ and one can define $\dg\f$ to be zero on the diagonal).

It is easy to see that if a function $\Phi$ on $\X\times\Y$ belongs to the {\it projective tensor
product}
$L^\be(E)\hat\otimes L^\be(F)$ of $L^\be(E)$ and $L^\be(F)$ (i.e., $\Phi$ admits a representation
$$
\Phi(\l,\mu)=\sum_{n\ge0}f_n(\l)g_n(\mu),
$$
where $f_n\in L^\be(E)$, $g_n\in L^\be(F)$, and
$$
\sum_{n\ge0}\|f_n\|_{L^\be}\|g_n\|_{L^\be}<\be),
$$
then $\Phi\in\fM(E,F)$, i.e., $\Phi$ is a Schur multiplier of the space of bounded linear operators. For such functions $\Phi$ we have
$$
\int\limits_\X\int\limits_\Y\Phi(\l,\mu)\,d E(\l)Q\,dF(\mu)=
\sum_{n\ge0}\left(\,\int\limits_\X f_n\,dE\right)Q\left(\,\int\limits_\Y g_n\,dF\right).
$$ 

More generally, $\Phi$ is a Schur multiplier  if $\Phi$ 
belongs to the {\it integral projective tensor product} $L^\be(E)\hat\otimes_{\rm i}
L^\be(F)$ of $L^\be(E)$ and $L^\be(F)$, i.e., $\Phi$ admits a representation
\bay
\label{ipt}
\Phi(\l,\mu)=\int_\O f(\l,\o)g(\mu,\o)\,d\s(\o),
\ey
where $(\O,\s)$ is a measure space, $f$ is a measurable function on $\X\times \O$,
$g$ is a measurable function on $\Y\times \O$, and
\bay
\label{ir}
\int_\O\|f(\cdot,\o)\|_{L^\be(E)}\|g(\cdot,\o)\|_{L^\be(F)}\,d\s(\o)<\be.
\ey
If $\Phi\in L^\be(E)\hat\otimes_{\rm i}L^\be(F)$, then
$$
\int\limits_\X\int\limits_\Y\Phi(\l,\mu)\,d E(\l)\,Q\,dF(\mu)=
\int\limits_\O\left(\,\int\limits_\X f(\l,\o)\,dE(\l)\right)Q
\left(\,\int\limits_\Y g(\mu,\o)\,dF(\mu)\right)\,d\s(\o).
$$
Clearly, the function 
$$
\o\mapsto \left(\,\int\limits_\X f(\l,\o)\,dE(\l)\right)Q
\left(\,\int\limits_\Y g(\mu,\o)\,dF(\mu)\right)
$$
is weakly measurable and
$$
\int\limits_\O\left\|\left(\,\int\limits_\X f(\l,\o)\,dE(\l)\right)T
\left(\,\int\limits_\Y g(\mu,\o)\,dF(\mu)\right)\right\|\,d\s(\o)<\be.
$$

Moreover, it can easily be seen that such functions $\Phi$ are Schur multipliers of any symmetrically normed ideal of operators.

It turns out that all Schur multipliers of the space of bounded linear operators can be obtained in this way. More precisely, the following result holds (see \cite{Pe1}):

\medskip

{\bf Theorem on Schur multipliers.} {\em Let $\Phi$ be a measurable function on 
$\X\times\Y$. The following are equivalent:

{\rm (i)} $\Phi\in\fM(E,F)$;

{\rm (ii)} $\Phi\in L^\be(E)\hat\otimes_{\rm i}L^\be(F)$;

{\rm (iii)} there exist measurable functions $f$ on $\X\times\O$ and $g$ on $\Y\times\O$ such that
{\em\rf{ipt}} holds and
\bay
\label{bs}
\left\|\int_\O|f(\cdot,\o)|^2\,d\s(\o)\right\|_{L^\be(E)}
\left\|\int_\O|g(\cdot,\o)|^2\,d\s(\o)\right\|_{L^\be(F)}<\be.
\ey
}

Note that the implication (iii)$\imp$(ii) was established in \cite{BS3}. Note also that
in the case of matrix Schur multipliers (this corresponds to discrete spectral measures
of multiplicity 1) the equivalence of (i) and (ii) was proved in \cite{Be}.

It is interesting to observe that if $f$ and $g$ satisfy \rf{ir}, then they also satisfy
\rf{bs}, but the converse is false. However, if $\Phi$ admits a representation of the form \rf{ipt}
with $f$ and $g$ satisfying \rf{bs}, then it also admits a (possibly different) representation of the
form \rf{ipt} with $f$ and $g$ satisfying \rf{ir}.

Note that in a similar way we can define the {\it projective tensor product} $A\hat\otimes B$
and the {\it integral projective tensor product} $A\hat\otimes_{\rm i} B$ of 
arbitrary Banach functions spaces $A$ and $B$.

Similar results also hold in the case of unitary operators.

\medskip

{\bf2.3. Sufficient conditions and necessary conditions.} We state here results in the case of unitary operators. Birman
and Solomyak proved in \cite{BS3} that if $\f$ is a function on the unit circle such that the derivative of $\f$ satisfies a H\"older condition of order $\a>0$ then the divided difference $\dg\f$ belongs to $\fM(E,F)$, which implies that if $U$ and $V$ are unitary operators, then
$$
\f(U)-\f(V)=\iint\limits_{\T\times\T}\frac{\f(\z)-\f(\t)}{\z-\t}\,dE_U(\z)\,(U-V)\,dE_V(\t),
$$
and so
$$
\|\f(U)-\f(V)\|\le\const\|U-V\|,
$$
i.e., $\f$ is {\it an operator Lipschitz function}. Moreover, it was shown in \cite{BS3} that under the same assumptions the function
$\f$ is {\it operator differentiable}, i.e., if $V=e^{{\rm i}A}U$, then the function $s\to\f\big(e^{{\rm i}sA}U\big)$ is differentiable
and
$$
\frac{d}{ds}\Big(\f(e^{{\rm i}sA}U)\Big)\Big|_{s=o}
={\rm i}\left(\iint\frac{\f(\z)-\f(\t)}{\z-\t}\,dE_U(\z)A\,dE_U(\t)\right)U.
$$

Later a much stronger result was obtained in \cite{Pe1}. It was shown in \cite{Pe1} that the same conclusions can be made under the assumption that $\f\in B_{\be1}^1$. Moreover, it was shown in \cite{Pe1} that if $\f\in B_{\be1}^1$, then $\f$ belongs to the projective tensor product $C(\T)\hat\otimes C(\T)$ of the space of continuous functions on $\T$ with itself, i.e., there exist functions
$f_n$ and $g_n$, $n\ge1$ in $C(\T)$ such that
$$
\sum_{n\ge1}\|f_n\|_\be\|g_n\|_\be<\be
$$
and
$$
\frac{\f(\z)-\f(\t)}{\z-\t}=\sum_{n\ge1}f_n(\z)g_n(\t).
$$
This implies that functions in $B_{\be1}^1$ are operator Lipschitz and operator differentiable. Moreover, this also implies that
functions in $B_{\be1}^1$ satisfy the inequality
$$
\|\f(U)-\f(V)\|_{\bS}\le\const\|U-V\|_\bS
$$
for arbitrary unitary operators $U$ and $V$ and for an arbitrary symmetrically normed ideal $\bS$.

Similar results hold for (not necessarily bounded) self-adjoint operators, see \cite{Pe4} and \cite{Pe8}.

It was proved in \cite{Pe1} that if $\f$ is an operator Lipschitz function on $\T$ (or if $\f$ is operator differentiable), then
$\f\in B_1^1$. This implies that
$$
\sum_{n\ge0}2^n\big|\hat\f\big(2^n\big)\big|<\be,
$$
and so the condition $\f\in C^1(\T)$ is not sufficient for $\f$ to be operator Lipschitz or operator differentiable.
An even stronger necessary condition found in \cite{Pe1} (see also \cite{Pe5} where that necessary condition was reformulated with the help of a remark by S. Semmes) says that if $\f$ is operator Lipschitz (or operator differentiable), then
$$
\f\in\cL\df\{\f:~|\n^2\f|\,dx\,dy~~\mbox{ is a Carleson measure in}~~\dd\},
$$
where $\n^2\f$ is the second gradient of
the harmonic extension of $\f$ to the unit disk. 

M. Frazier observed that actually $\cL$ is the Triebel--Lizorkin space $F^1_{\be1}$. Note that the definition of the
Triebel--Lizorkin spaces $\dot F^s_{pq}$ on $\R^n$ for $p=\be$ and  $q>1$ can be found in \cite{T}, \S\,5.1.
A definition for all $q>0$, which is equivalent to Triebel's
definition when $q>1$, was given by Frazier and Jawerth in \cite{FrJ}. 
Their approach did not use harmonic extensions, but a
straightforward exercise in comparing kernels shows that Frazier and
Jawerth's definition of $\dot F^1_{\infty1}$ is equivalent to the
definition requiring $|\n^2 u| dx dy$ to be a Carleson measure on the
upper half-space.  Our space $\cL$ is the analogue for the unit disc.

It was observed in \cite{Pe2} that if $\f$ is an analytic function in $\big(B_{\be1}^1\big)_+$, then the divided difference
$\dg\f$ belongs to the projective tensor product $C_A\hat\otimes C_A$ of the disk-algebra $C_A$ with itself, and so if 
$\f\in\big(B_{\be1}^1\big)_+$, then for arbitrary contractions $T$ and $R$ and an arbitrary symmetrically normed ideal $\bS$
the following inequality holds:
$$
\|\f(T)-\f(R)\|_\bS\le\const\|T-R\|_\bS.
$$

Recently in \cite{KS} it was proved that if $\f\in C_A$, the the following are equivalent:

(i) $\|\f(U)-\f(V)\|\le\const\|U-V\|$ for arbitrary unitary operators $U$ and $V$;

(ii) $\|\f(T)-\f(R)\|\le\const\|T-R\|$ for arbitrary contractions $T$ and $R$.

\medskip

{\bf 2.4. Multiple operator integrals.} The equivalence of (i) and (ii) in the Theorem on Schur multipliers suggests the idea 
explored in  \cite{Pe8} to define multiple operator integrals. 

To simplify the notation, we consider here the case of triple operator integrals; the case of arbitrary multiple operator integrals can be treated in the same way.

Let $(\X,E)$, $(\Y,F)$, and $(\cZ,G)$
be spaces with spectral measures $E$, $F$, and $G$ on Hilbert spaces $\h_1$, $\h_2$, and $\h_3$. Suppose that
$\psi$ belongs to the integral projective tensor product
$L^\be(E)\hat\otimes_{\rm i}L^\be(F)\hat\otimes_{\rm i}L^\be(G)$, i.e., $\psi$ admits a representation
\bay
\label{ttp}
\Phi(\l,\mu,\nu)=\int_\O f(\l,\o)g(\mu,\o)h(\nu,\o)\,d\s(\o),
\ey
where $(\O,\s)$ is a measure space, $f$ is a measurable function on $\X\times \O$,
$g$ is a measurable function on $\Y\times \O$, $h$ is a measurable function on $\cZ\times \O$,
and
$$
\int_Q\|f(\cdot,x)\|_{L^\be(E)}\|g(\cdot,x)\|_{L^\be(F)}\|h(\cdot,x)\|_{L^\be(G)}\,d\s(x)<\be.
$$

Suppose now that $T_1$ is a bounded linear operator from $\h_2$ to $\h_1$ and $T_2$ is a bounded linear operator from $\h_3$ to $\h_2$. For a function $\Phi$ in
$L^\be(E)\hat\otimes_{\rm i}L^\be(F)\hat\otimes_{\rm i}L^\be(G)$ of the form \rf{ttp}, we put
\begin{align}
\label{opr}
&\int\limits_\X\int\limits_\Y\int\limits_\cZ\Phi(\l,\mu,\nu)
\,d E(\l)T_1\,dF(\mu)T_2\,dG(\nu)\nonumber\\[.2cm]
\df&\int\limits_\O\left(\,\int\limits_\X f(\l,\o)\,dE(\l)\right)T_1
\left(\,\int\limits_\Y g(\mu,\o)\,dF(\mu)\right)T_2
\left(\,\int\limits_\cZ h(\nu,\o)\,dG(\nu)\right)\,d\s(\o).
\end{align}

The following lemma from \cite{Pe8} (see also \cite{ACDS} for a different proof) shows that the definition does not depend on the choice of a representation \rf{ttp}.

\begin{lem}
\label{kor}
Suppose that $\psi\in L^\be(E)\hat\otimes_{\rm i}L^\be(F)\hat\otimes_{\rm i}L^\be(G)$. Then the 
right-hand side of {\em\rf{opr}} does not depend on the choice of a representation {\em\rf{ttp}}.
\end{lem}

It is easy to see that the following inequality holds
$$
\left\|\int\limits_\X\int\limits_\Y\int\limits_\cZ\psi(\l,\mu,\nu)
\,dE(\l)T_1\,dF(\mu)T_2\,dG(\nu)\right\|
\le\|\psi\|_{L^\be\hat\otimes_{\rm i}L^\be\hat\otimes_{\rm i}L^\be}\cdot\|T_1\|\cdot\|T_2\|.
$$

In particular, the triple operator integral on the left-hand side of \rf{opr} can be defined if $\Phi$ belongs to the projective
tensor product $L^\be(E)\hat\otimes L^\be(F)\hat\otimes L^\be(G)$, i.e. $\Phi$ admits a representation
$$
\Phi(\l,\mu,\nu)=\sum_{n\ge1}f_n(\l)g_n(\mu)h_n(\nu),
$$
where $f_n\in L^\be(E)$, $g_n\in L^\be(F)$, $h_n\in L^\be(G)$ and
$$
\sum_{n\ge1}\|f_n\|_{L^\be(E)}\|g_n\|_{L^\be(F)}|\|h_n\|_{L^\be(G)}<\be.
$$

In a similar way one can define multiple operator integrals, see \cite{Pe8}.

For a function $\f$ on the circle the  {\it divided differences $\dg^k\f$ of order $k$} are defined inductively as follows:
$$
\dg^0\f\df\f;
$$
if $k\ge1$, then in the case $\l_1,\l_2,\cdots,\l_{k+1}$ are distinct points in $\T$,
$$
(\dg^{k}\f)(\l_1,\cdots,\l_{k+1})\df
\frac{(\dg^{k-1}\f)(\l_1,\cdots,\l_{k-1},\l_k)-
(\dg^{k-1}\f)(\l_1,\cdots,\l_{k-1},\l_{k+1})}{\l_{k}-\l_{k+1}}
$$
(the definition does not depend on the order of the variables). Clearly,
$$
\dg\f=\dg^1\f.
$$
If $\f\in C^k(\T)$, then $\dg^{k}\f$ extends by continuity to a function defined for all points $\l_1,\l_2,\cdots,\l_{k+1}$.

It was shown in \cite{Pe8} that if $\f\in B_{\be1}^n$, then $\dg^{n}\f$ belongs to the projective tensor product
$\underbrace{C(\T)\hat\otimes \cdots\hat\otimes C(\T)}_{n+1}$. Moreover, it was shown in \cite{Pe8}
if $U$ and $V$ are unitary operators, $V=e^{iA}U$, then the function
$$
t\mapsto e^{{\rm i}tA}U
$$
has $n$ derivatives in the norm and
\begin{align*}
&\frac{d^n}{dt^n}\Big(\big(\f(e^{{\rm i}tA}U\big)\Big)\Big|_{s=0}\\[.2cm]
=&{\rm i}^nn!\left(\underbrace{\int\cdots\int}_{n+1}(\dg^n\f)(\l_1,\cdots,\l_{n+1})
\,dE_U(\l_1)A\cdots A\,dE_U(\l_{n+1})\right)U^n.
\end{align*}

The reasoning given in \cite{Pe8} shows that 
\bay
\f\in\big(B_{\be1}^n\big)_+\quad\imp\quad \dg^{n}\f \in\underbrace{C_A\hat\otimes \cdots\hat\otimes C_A}_{n+1}
\ey

Note that recently in \cite{JTT} Haagerup tensor products were used to define multiple operator integrals. However, it is
not clear whether this can lead to a broader class of functions $\f$, for which $\dg^n\f$ can be integrated.

\medskip

{\bf2.5.~Semi-spectral measures.} Let $\h$ be a Hilbert space and let $(\X,\B)$ be a measurable space.
A map $\E$ from $\B$ to the algebra $B(\h)$ of all bounded operators on $\h$ is called a {\it semi-spectral measure}
if 
$$
\E(\D)\ge\0,\quad\D\in\B,
$$
$$
\E(\varnothing)=\0\quad\mbox{and}\quad\E(\X)=I,
$$
and for a sequence $\{\D_j\}_{j\ge1}$ of disjoint sets in $\B$,
$$
\E\left(\bigcup_{j=1}^\be\D_j\right)=\lim_{N\to\be}\sum_{j=1}^N\E(\D_j)\quad\mbox{in the weak operator topology}.
$$

If $\K$ is a Hilbert space, $(\X,\B)$ is a measurable space,  $E:\B\to B(\K)$ is a spectral measure, and $\h$ is
a subspace of $\K$, then it is easy to see that the map $\E:\B\to B(\h)$ defined by
\bay
\label{dil}
\E(\D)=P_\h E(\D)\big|\h,\quad\D\in\B,
\ey
is a semi-spectral measure. Here $P_\h$ stands for the orthogonal projection onto $\h$.

Naimark proved in \cite{N}  (see also \cite{SNF}) that all semi-spectral measures can be obtained in this way, i.e.,
a semi-spectral measure is always a {\it compression} of a spectral measure. A spectral measure $E$ satisfying \rf{dil} is called a {\it spectral dilation of the semi-spectral measure} $\E$.

A spectral dilation $E$ of a semi-spectral measure $\E$ is called {\it minimal} if 
$$
\K=\clos\spn\{E(\D)\h:~\D\in\B\}.
$$

It was shown in \cite{MM} that if $E$ is a minimal spectral dilation of a semi-spectral measure $\E$, then
$E$ and $\E$ are mutually absolutely continuous and all minimal spectral dilations of a semi-spectral measure are isomorphic in the natural sense.

If $\f$ is a bounded complexed-valued measurable function on $\X$ and $\E:\B\to B(\h)$ is a semi-spectral measure, then the integral
\bay
\label{iss}
\int_\X f(x)\,d\E(x)
\ey
can be defined as
\bay
\label{voi}
\int_\X f(x)\,d\E(x)=\left.P_\h\left(\int_\X f(x)\,d E(x)\right)\right|\h,
\ey
where $E$ is a spectral dilation of $\E$. It is easy to see that the right-hand side of \rf{voi} does not depend on the choice
of a spectral dilation. The integral \rf{iss} can also be computed as the limit of sums
$$
\sum f(x_\a)\E(\D_\a),\quad x_\a\in\D_\a,
$$
over all finite measurable partitions $\{\D_\a\}_\a$ of $\X$.

If $T$ is a contraction on a Hilbert space $\h$, then by the Sz.-Nagy dilation theorem
(see \cite{SNF}),  $T$ has a unitary dilation, i.e., there exist a Hilbert space $\K$ such that
$\h\subset\K$ and a unitary operator $U$ on $\K$ such that
\bay
\label{DT}
T^n=P_\h U^n\big|\h,\quad n\ge0,
\ey
where $P_\h$ is the orthogonal projection onto $\h$. Let $E_U$ be the spectral measure of $U$.
Consider the operator set function $\E$ defined on the Borel subsets of the unit circle $\T$ by
$$
\E(\D)=P_\h E_U(\D)\big|\h,\quad\D\subset\T.
$$
Then $\E$ is a semi-spectral measure. It follows immediately from
\rf{DT} that 
\bay
\label{step}
T^n=\int_\T \z^n\,d\E(\z)=P_\h\int_\T\z^n\,dE_U(\z)\Big|\h,\quad n\ge0.
\ey
Such a semi-spectral measure $\E$ is called a {\it semi-spectral measure} of $\T$.
Note that it is not unique. To have uniqueness, we can consider a minimal unitary dilation $U$ of $T$,
which is unique up to an isomorphism (see \cite{SNF}).

It follows easily from from \rf{step} that 
$$
\f(T)=P_\h\int_\T\f(\z)\,dE_U(\z)\Big|\h
$$
for an arbitrary function $\f$ in the disk-algebra $C_A$.

\

\section{\bf Double operator integrals with respect to semi-spectral measures}
\setcounter{equation}{0}

\

In this section we extend the Birman--Solomyak theory of double operator integrals to the case of semi-spectral measures.

Suppose that $(\X_1,\B_1)$ and $(\X_2,\B_2)$ are measurable spaces, $\h_1$ and $\h_2$ are Hilbert spaces, and
$\E_1:\B_1\to B(\h_1)$ and $\E_2:\B_2\to B(\h_2)$ are semi-spectral measures.
For a bounded measurable function $\Phi$ on $\X_1\times\X_2$ and an operator $Q:\h_2\to\h_1$ we consider double operator integrals
\bay
\label{dois}
\iint\limits_{\X_1\times\X_2}\Phi(x_1,x_2)\,d\E_1(x_1)Q\,d\E_2(X_2).
\ey
In the case when $Q$ is a Hilbert--Schmidt operator, integrals of the form \rf{dois} can be interpreted as in the case of double operator integrals with respect to spectral measures (see \S\,2.2). Indeed, we define the map $\F$ on the set of all measurable rectangles
$\D_1\times\D_2$ by
$$
\F(\D_1\times\D_2)Q=\E_1(\D_1)Q\E_2(\D_2),\quad Q\in\bS_2(\h_2,\h_1).
$$
Clearly, $\F(\D_1\times\D_2)$ is a bounded linear operator on $\bS_2(\h_2,\h_1)$ that satisfies the inequalities $\0\le\F(\D_1\times\D_2)\le I$. 

\begin{lem}
\label{ext}
$\F$ extends to a semi-spectral measure on $\B_1\times\B_2$.
\end{lem}

\Pf  Let $E_1$ and $E_2$ be spectral dilations of $\E_1$ and $\E_2$ on Hilbert spaces $\K_1$ and $\K_2$. Define the map $F$ on measurable rectangles $\D_1\times\D_2$ by
$$
F(\D_1\times\D_2)Q=E_1(\D_1)QE_2(\D_2),\quad Q\in\bS_2(\h_2,\h_1).
$$
By the theorem of Birman and Solomyak (see \cite{BS4}), $F$ extends to a spectral measure (which will also be denoted by $F$) defined on $\X_1\times\X_2$.
Clearly,
$$
\F(\D_1\times\D_2)Q=P_{\h_1}F(\D_1\times\D_2)\big(P_{\h_1}QP_{\h_2}\big)\Big|\h_2,\quad Q\in\bS_2(\h_2,\h_1).
$$
We can define now the operator $\F(\O)$ for an arbitrary $\O\in\X_1\times\X_2$ by
$$
\F(\O)Q=P_{\h_1}F(\O)\big(P_{\h_1}QP_{\h_2}\big)\Big|\h_2,\quad Q\in\bS_2(\h_2,\h_1).
$$
It is easy to see that $\F$ is a semi-spectral measure on the Hilbert space $\bS_2(\h_2,\h_1)$. $\bl$

Given a bounded measurable function $\Phi$ on $\X_1\times\X_2$ and a Hilbert--Schmidt operator $Q:\h_2\to\h_1$,
we can define the double operator integral
$$
\iint\limits_{\X_1\times\X_2}\Phi(x_1,x_2)\,d\E_1(x_1)Q\,d\E_2(x_2)
$$
as in the case of integration with respect to spectral measures:
$$
\iint\limits_{\X_1\times\X_2}\Phi(x_1,x_2)\,d\E_1(x_1)Q\,d\E_2(x_2)
\df\left(\int_{\X_1\times\X_2}\Phi\,d\F\right)Q.
$$
It is easy to see that 
$$
\iint\limits_{\X_1\times\X_2}\Phi(x_1,x_2)\,d\E_1(x_1)Q\,d\E_2(x_2)=
P_{\h_1}\iint\limits_{\X_1\times\X_2}\Phi(x_1,x_2)\,dE_1(x_1)Q\,dE_2(x_2)\big|\h_2.
$$
Clearly, 
\bay
\label{fn2}
\left\|~\,\iint\limits_{\X_1\times\X_2}\Phi(x_1,x_2)\,d\E_1(x_1)Q\,d\E_2(x_2)\right\|_{\bS_2}
\le\sup_{x_1\in\X_1,x_2\in\X_2}
|\Phi(x_1,x_2)|\cdot\|Q\|_{\bS_2}.
\ey

We need the following fact.

\begin{lem}
\label{poss}
Suppose that $\E_1$ and $\E_2$ are semi-spectral measure as above and $Q$ is a Hilbert-Schmidt operator.
If $\{\Phi_n\}_{n\ge1}$ is a sequence of measurable functions such that 
$$
\lim_{n\to\be}\Phi_n(\l,\mu)=\Phi(\l,\mu),\quad \l\in\X_1,~\mu\in\X_2,
$$
and
$$
\sup_n\sup_{\l,\,\mu\in\X_1\times\X_2}|\Phi_n(\l,\mu)|<\be,
$$
then
$$
\lim_{n\to\be}\left\|~\int\limits_{\X_1}\int\limits_{\X_2}\Phi_n(\l,\mu)\,d E_1(\l)\,Q\,dE_2(\mu)
-\int\limits_{\X_1}\int\limits_{\X_2}\Phi(\l,\mu)\,d E_1(\l)\,Q\,dE_2(\mu)\right\|_{\bS_2}=0.
$$
\end{lem}

\Pf The result follows immediately from the same fact in the case of spectral measures, see \rf{poto}. $\bl$

Let us proceed now to double operator integrals \rf{dois} with bounded operators $Q$.

Suppose that $\E_1$ and $\E_2$ are semi-spectral measures and $E_1$ and $E_2$ are their minimal spectral dilations.
If $\Phi$ is a Schur multiplier of the space of bounded linear operators, then the double operator integral
\bay
\label{doss}
\iint\limits_{\X_1\times\X_2}\Phi(x_1,x_2)\,d\E_1(x_1)Q\,d\E_2(x_2)
\ey
is defined as 
$$
P_{\h_1}\iint\limits_{\X_1\times\X_2}\Phi(x_1,x_2)\,dE_1(x_1)Q\,dE_2(x_2)\big|\h_2.
$$


It is easy to see that if $\Phi\in L^\be(\E_1)\hat\otimes_{\rm i}L^\be(\E_2)$
and
$$
\Phi(x_1,x_2)=\int_\O f(x_1,\o)g(x_2,\o)\,d\s(\o)
$$
with
$$
\int_\O\|f(\cdot,\o)\|_{L^\be(\E_1)}\|g(\cdot,\o)\|_{L^\be(\E_2)}\,d\s(\o)<\be,
$$
then $\Phi$ is a Schur multiplier of the space of bounded linear operators, the integral \rf{doss} is equal to
$$
\int\limits_\O\left(\,\int\limits_{\X_1}f(x_1,\o)\,d\E_1(x_1)\right)Q
\left(\,\int\limits_{\X_2}\!g(x_2,\o)\,d\E_2(x_2)\right)\,d\s(\o)
$$
and its norm is less than or equal to
$$
\|Q\|\int_\O\|f(\cdot,\o)\|_{L^\be(\E_1)}\|g(\cdot,\o)\|_{L^\be(\E_2)}\,d\s(\o).
$$

We can define now multiple operator integrals in the same way as in the case of semi-spectral measures.
Suppose that $\h_1,\h_2,\cdots,\h_n$ are a Hilbert spaces and for \lb$j=1,2,\cdots,n-1$, $Q_j$ 
is a bounded linear operator from $\h_{j+1}$ to $\h_j$. Suppose also that
 $\E_1,\cdots,\E_n$ 
are semi-spectral measures defined on $\s$-algebras of $\X_1,\cdots,\X_n$, $\E_j$ takes valued in the space $\B(\h_j)$,
and $\Phi$ is a function on $\X_1\times\cdots\times\X_n$ of class 
$L^\be(\E_1)\hat\otimes_{\rm i}\cdots\hat\otimes_{\rm i}L^\be(\E_n)$, i.e., $\Phi$ admits a representation
\bay
\label{itr}
\Phi(x_1,\cdots,x_n)=\int_\O f_1(x_1,\o)\cdots f_n(x_n,\o)\,d\s(\o),
\ey
in which
\bay
\label{itco}
\int_\O\|f_1(,\cdot,\o)\|_{L^\be(\E_1)}\cdots\|f_n(,\cdot,\o)\|_{L^\be(\E_n)}\,d\s(\o)<\be.
\ey
We define the multiple operator integral
$$
\underbrace{\int\cdots\int}_n\Phi(x_1,\cdots,x_n)
\,d\E_1(x_1)Q_1\cdots Q_{n-1}\,d\E_n(x_n)
$$
as
\bay
\label{mois}
\int_\O\left(\int_{\X_1}f_1(x_1,\o)\,d\E_1(x_1)\right)Q_1\cdots Q_{n-1}\left(\int_{\X_n}f_1(x_n,\o)\,d\E_n(x_n)\right)\,d\s(\o).
\ey
Certainly, we have to prove that the multiple operator integral is well defined. In other words, we have to show that the value of
\rf{mois} does not depend on the choice of a representation \rf{itr}, which is a consequence of the following lemma.

\begin{lem}
\label{nez}
Suppose that 
$$
\int_\O f_1(x_1,\o)\cdots f_n(x_n,\o)\,d\s(\o)=0,\quad x_1\in\X_1,\cdots,x_n\in\X_n,
$$
and {\em\rf{itco}} holds. Then
\bay
\label{nol}
\int_\O\left(\int_{\X_1}f_1(x_1,\o)\,d\E_1(x_1)\right)Q_1\cdots Q_{n-1}\left(\int_{\X_n}f_n(x_n,\o)\,d\E_n(x_n)\right)\,d\s(\o)
\ey
is the zero operator.
\end{lem}

\Pf We deduce the lemma from the corresponding fact for multiple operator integrals with respect to spectral measures.
Suppose that $E_j$ is a minimal spectral dilation of $\E_j$ and $E_j$ takes values in $\B(\K)$. Then the integral in \rf{nol} is equal to
\bay
\label{korf}
P_\h\left(
\int_\O\left(\int_{\X_1}\!f_1(x_1,\o)\,dE_1(x_1)\right)\!B_1\cdots B_{n-1}\!\left(\int_{\X_n}\!f_n(x_n,\o)\,dE_n(x_n)\right)\,d\s(\o)
\right)\Big|\h,
\ey
where $B_j\df P_{\h_j} Q_jP_{\h_{j+1}}$, $1\le j\le n-1$. It follows now from Lemma \ref{kor} that the operator in \rf{korf} is the zero operator. $\bl$

\

\section{\bf Analogs of  the Birman--Solomyak formulae}
\setcounter{equation}{0}

\

In this section we obtain analogs of the Birman--Solomyak formulae \rf{BSdof} and \rf{BScom} for contractions. 

If $T$ is a contraction on Hilbert space, then $\E_T$ denotes a semi-spectral measure of $\T$.

Recall that if $\f'\in C_A$, the function $\cd\f$ extends to the diagonal
$$
\bs{\D}\df\big\{(\z,\z):~\z\in\T\big\}
$$
by continuity: $(\cd\f)(\z,\z)=\f'(\z)$, $\z\in\T$.

\begin{thm}
\label{BST}
Let $\f\in\left(B_{\be1}^1\right)_+$. Then for contractions $T$ and $R$ on Hilbert space the following formula holds:
\bay
\label{BSF}
\f(T)-\f(R)=\iint\limits_{\T\times\T}\frac{\f(\z)-\f(\t)}{\z-\t}\,d\E_T(\z)\,(T-R)\,d\E_R(\t).
\ey
\end{thm}

Recall that it follows from the results of \cite{Pe1} that the function $\cd\f$ belongs to the projective tensor product $C_A\hat{\otimes}C_A$
(see \S\,2.2),
and so the right-hand side of \rf{BSF} is well defined.

\medskip

{\bf Proof of Theorem \ref{BST}.}  Suppose that
$$
(\cd\f)(\z,\t)=\frac{\f(\z)-\f(\t)}{\z-\t}=\sum_{n\ge0}f_n(\z)g_n(\t),
$$
where $f_n\in C_A$, $g_n\in C_A$, and
$$
\sum_{n\ge0}\|f_n\|_\be\|g_n\|_\be<\be.
$$
We have
\begin{align*}
\iint\limits_{\T\times\T}\frac{\f(\z)-\f(\t)}{\z-\t}\,d\E_T(\z)(T-R)\,d\E_R(\t)&=
\sum_{n\ge0}f_n(T)(T-R)g_n(R)\\[.2cm]
&=\sum_{n\ge0}Tf_n(T)g_n(R)-\sum_{n\ge0}f_n(T)g_n(R)R.
\end{align*}
Clearly, 
$$
\sum_{n\ge0}Tf_n(T)g_n(R)=\iint\limits_{\T\times\T}\z f_n(\z)g_n(\t)\,d\E_T(\z)\,d\E_R(\t)
=\iint\limits_{\T\times\T}\z(\cd\f)(\z,\t)\,d\E_T(\z)\,d\E_R(\t)
$$
and
$$
\sum_{n\ge0}f_n(T)g_n(R)R=\iint\limits_{\T\times\T} f_n(\z)g_n(\t)\t\,d\E_T(\z)\,d\E_R(\t)
=\iint\limits_{\T\times\T} \t(\cd\f)(\z,\t)\,d\E_T(\z)\,d\E_R(\t).
$$
Thus
\begin{align*}
\iint\limits_{\T\times\T}(\cd\f)(\z,\t)\,d\E_T(\z)(T-R)\,d\E_R(\t)&=
\iint\limits_{\T\times\T}\z(\cd\f)(\z,\t)\,d\E_T(\z)\,d\E_R(\t)\\[.2cm]
&-
\iint\limits_{\T\times\T} \t(\cd\f)(\z,\t)\,d\E_T(\z)\,d\E_R(\t)\\[.2cm]
&=\iint\limits_{\T\times\T} (\z-\t)(\cd\f)(\z,\t)\,d\E_T(\z)\,d\E_R(\t)\\[.2cm]
&=\iint\limits_{\T\times\T}\big(\f(\z)-\f(\t)\big)\,d\E_T(\z)\,d\E_R(\t)\\[.2cm]
&=\f(T)-\f(R).\quad\bl
\end{align*}

Let us consider the case when $T-R\in\bS_2$. The following result establishes formula \rf{BSF} for functions with derivatives
in $C_A$. 

\begin{thm}
\label{S2}
Let $\f$ be a function analytic in $\dd$ such that $\f'\in C_A$. 
If $T$ and $R$ are contractions such that $T-R\in\bS_2$, then formula
{\em\rf{BSF}} holds.
\end{thm}

\Pf Let $\{\f_n\}_{n\ge1}$ be a sequence of polynomials such that
$$
\|\f'-\f'_n\|_\be\to0\quad\mbox{as}\quad n\to\be.
$$
It follows from Theorem \ref{BST} that 
$$
\f_n(T)-\f_n(R)=\iint\limits_{\T\times\T}\frac{\f_n(\z)-\f_n(\t)}{\z-\t}\,d\E_T(\z)\,(T-R)\,d\E_R(\t).
$$
By von Neumann's inequality, 
$$
\lim_{n\to\be}\|\f_n(T)-\f(T)\|\to0\quad\mbox{and}\quad\lim_{n\to\be}\|\f_n(R)-\f(R)\|\to0.
$$
The result follows now from the trivial observation
that
$$
\sup_{\z,\t\in\T}|(\cd\f_n)(\z,\t)-(\cd\f)(\z,\t)|\to0\quad\mbox{as}\quad n\to\be
$$
which implies that
$$
\lim_{n\to\be}\left\|\,\,\iint\limits_{\T\times\T}\Big((\cd\f_n)(\z,\t)-(\cd\f)(\z,\t)\Big)\,d\E_T(\z)\,(T-R)\,d\E_R(\t)
\right\|_{\bS_2}=0.\quad\bl
$$

Consider now the more general case when $\f'\in H^\be$.

\begin{thm}
\label{mgt}
Let $\f'\in H^\be$ and let $\Phi$ be the function on $\T\times\T$ defined by
\bay
\label{Phi}
\Phi\Big|\big(\T\times\T\setminus\bs{\D}\big)=(\cd\f)\Big|\big(\T\times\T\setminus\bs{\D}\big)\quad
\mbox{and}\quad\Phi\Big|\bs{\D}=\0.
\ey
If $T$ and $R$ are contractions such that $T-R\in\bS_2$, then
$$
\f(T)-\f(R)=\iint\limits_{\T\times\T}\Phi(\z,\t)\,d\E_T(\z)\,(T-R)\,d\E_R(\t).
$$
\end{thm}

\Pf Let $\f_n$ be the $n$th Ces\'aro mean of the Fourier series of $\f$.
Then $\f_n\in\big(B_{\be1}^1\big)_+$ and by Theorem \ref{BST},
\begin{align*}
\f_n(T)-\f_n(R)&=\iint\limits_{\T\times\T}\frac{\f_n(\z)-\f_n(\t)}{\z-\t}\,d\E_T(\z)\,(T-R)\,d\E_R(\t)\\[.2cm]
&=
\iint\limits_{\T\times\T}\big(\f_n(\z)-\f_n(\t)\big)\,d\E_T(\z)\,d\E_R(\t)\\[.2cm]
&=\iint\limits_{(\T\times\T)\setminus\bs{\D}}\big(\f_n(\z)-\f_n(\t)\big)\,d\E_T(\z)\,d\E_R(\t).
\end{align*}

The same reasoning as in the proof of Theorem \ref{BST} shows that 
$$
\iint\limits_{(\T\times\T)\setminus\bs{\D}}\frac{\f_n(\z)-\f_n(\t)}{\z-\t}\,d\E_T(\z)\,(T-R)\,d\E_R(\t)=
\iint\limits_{(\T\times\T)\setminus\bs{\D}}\big(\f_n(\z)-\f_n(\t)\big)\,d\E_T(\z)\,d\E_R(\t).
$$

Clearly, $\|\f_n-\f\|_\be\to0$ as $n\to\be$, and so by von Neumann's inequality,
$$
\f_n(T)-\f_n(R)\to\f(T)-\f(R)
$$
in the operator norm.

On the other hand, since
$$
\lim_{n\to\be}\frac{\f_n(\z)-\f_n(\t)}{\z-\t}=\frac{\f(\z)-\f(\t)}{\z-\t},\quad\z,\,\t\in\T,
$$
and
$$
\sup_n\sup_{\z,\t\in\T}\left|\frac{\f_n(\z)-\f_n(\t)}{\z-\t}\right|<\be,
$$
it follows from Lemma \ref{poss} that
$$
\left\|~\,\,\iint\limits_{(\T\times\T)\setminus\bs{\D}}\left(\frac{\f_n(\z)-\f_n(\t)}{\z-\t}-\frac{\f(\z)-\f(\t)}{\z-\t}\right)\,d\E_T(\z)\,(T-R)\,d\E_R(\t)
\right\|_{\bS_2}\to0
$$
as $n\to\be$.

To complete the proof, it remains to observe that
$$
\iint\limits_{(\T\times\T)\setminus\bs{\D}}\frac{\f(\z)-\f(\t)}{\z-\t}\,d\E_T(\z)\,(T-R)\,d\E_R(\t)=
\iint\limits_{\T\times\T}\Phi(\z,\t)\,d\E_T(\z)\,(T-R)\,d\E_R(\t).\quad\bl
$$

The following result is an immediate consequence of Theorem \ref{mgt}; it was obtained recently in \cite{KS} by a
completely different method.

\begin{cor}
\label{hsp}
Suppose that $\f$ is a function analytic in $\dd$ such that $\f'\in H^\be$. If $T$ and $R$ are contractions on
Hilbert space such that $T-R\in\bS_2$, then
$$
\f(T)-\f(R)\in\bS_2\quad\mbox{and}\quad\|\f(R)-\f(T)\|_{\bS_2}\le\|\f'\|_{H^\be}\|T-R\|_{\bS_2}.
$$
\end{cor}

Let us obtain now analogs of formula \rf{BScom}.

\begin{thm}
\label{BS2Bes}
Let $\f\in\big(B^1_{\be1}\big)_+$. Then for a contraction $T$ and a bounded linear operator $Q$ on Hilbert space the following formula holds:
$$
\f(T)Q-Q\f(T)=\iint\limits_{\T\times\T}\frac{\f(\z)-\f(\t)}{\z-\t}\,d\E_T(\z)\,(TQ-QT)\,d\E_T(\t).
$$
\end{thm}

The following result is a Hilbert--Schmidt version of Theorem \ref{BS2Bes}.

\begin{thm}
\label{S22}
Let $\f'\in H^\be$ and let $\Phi$ be the function on $\T\times\T$ defined by {\em\rf{Phi}}. Suppose that $\T$ is a contraction and $Q$ is a bounded linear operator such that
\bay
\label{com2}
TQ-QT\in\bS_2.
\ey
Then 
$$
\f(T)Q-Q\f(T)=\iint\limits_{\T\times\T}\Phi(\z,\t)\,d\E_T(\z)\,(TQ-QT)\,d\E_T(\t).
$$
\end{thm}

Theorems \ref{BS2Bes} and \ref{S22} can be proved in the same way as Theorems \ref{BST} and \ref{mgt}.

The following inequality is an immediate consequence of Theorem \ref{S22}; recently it was proved in \cite{KS} by a different method.

\begin{cor}
\label{sled}
Suppose that $\f'\in H^\be$, $T$ is a contraction and $Q$ is a bounded linear operator on Hilbert space such that
{\em\rf{com2}} holds. Then
$$
\f(T)Q-Q\f(T)\in\bS_2
$$
and
$$
\|\f(T)Q-Q\f(T)\|_{\bS_2}\le\|\f'\|_{H^\be}\|TQ-QT\|_{\bS_2}.
$$
\end{cor}

\

\section{\bf Differentiability of operator functions in the operator norm}
\setcounter{equation}{0}

\

Let $T$ and $R$ be contractions on a Hilbert space $\h$. For $t\in[0,1]$, consider the operator $T_t$
defined by
\bay 
\label{Tt}
T_t=(1-t)T+tR.
\ey
Clearly, $T_t$ is a contraction. 

In this section for a functions $\f\in C_A$ we consider the problem of differentiability of 
the function
\bay
\label{fun}
t\mapsto \f(T_t),\quad0\le t\le1.
\ey
and the problem of the existence of higher derivatives. We also compute the derivatives in terms of multiple operator
integrals with respect to semi-spectral measures.
 
Let $\E_t$ be a semi-spectral measure of $T_t$ on the unit circle $\T$,
i.e.,
$$
T_t^n=\int\limits_\T\z^n\,d\E_t(\z),\quad n\ge0.
$$
Put $\E\df\E_0$.

\begin{thm}
\label{per}
Suppose that $\f\in\big(B^1_{\be1}\big)_+$. Then the function {\em\rf{fun}} is differentiable 
in the norm and
$$
\frac{d}{ds}\f(T_s)\Big|_{s=t}=\iint\limits_{\T\times\T}\frac{\f(\z)-\f(\t)}{\z-\t}\,d\E_t(\z)\,(R-T)\,d\E_t(\t).
$$
\end{thm}

\Pf We have 
\begin{align*}
\frac1s\big(\f(T_{t+s})-\f(T_t)\big)&=
\frac1s\iint\limits_{\T\times\T}\frac{\f(\z)-\f(\t)}{\z-\t}\,d\E_{t+s}(\z)\,(T_{t+s}-T_t)\,d\E_t(\t)\\[.2cm]
&=\iint\limits_{\T\times\T}\frac{\f(\z)-\f(\t)}{\z-\t}\,d\E_{t+s}(\z)\,(R-T)\,d\E_t(\t).
\end{align*}

As we have mentioned in \S\,2.3, it was shown in \cite{Pe1} that  $\dg\f$ admits a representation
\bay
\label{ten}
\frac{\f(\z)-\f(\t)}{\z-\t}=\sum_nf_n(\z)g_n(\t),
\ey
where $f_n,\,g_n\in C_A$ and
\bay
\label{pro}
\sum_n\|f_n\|_\be\|g_n\|_\be<\be.
\ey
The following identities hold:
\begin{align*}
\iint\limits_{\T\times\T}\frac{\f(\z)-\f(\t)}{\z-\t}\,d\E_{t+s}(\z)\,(R-T)\,d\E_t(\t)&=
\sum_n\iint \limits_{\T\times\T}f_n(\z)g_n(\t)\,d\E_{t+s}(\z)\,(R-T)\,d\E_t(\t)\\[.2cm]
&=\sum_n f_n(T_{t+s})(R-T)g_n(T_t).
\end{align*}
Clearly,
$$
\lim_{s\to0}\|f_n(T_{t+s})-f_n(T)\|=0
$$
and in view of \rf{pro},
$$
\lim_{s\to0}\frac1s\big(\f(T_{t+s})-\f(T_t)\big)=\sum_n f_n(T_{t})(R-T)g_n(T_t)
$$
in the norm. It remains to observe that by \rf{ten},
\begin{align*}
\lim_{s\to0}
\sum_n f_n(T_{t+s})(R-T)g_n(T_t)&=\sum_n f_n(T_{t})(R-T)g_n(T_t)\\[.2cm]
&=\sum_n\iint\limits_{\T\times\T} f_n(\z)g_n(\t)\,d\E_t(\z)\,(R-T)\,d\E_t(\t)\\[.2cm]
&=\iint\limits_{\T\times\T}\frac{\f(\z)-\f(\t)}{\z-\t}\,d\E_t(\z)\,(R-T)\,d\E_t(\t).\quad\bl
\end{align*}

In what follows to simplify the notation, we will not specify the domain of integration: all double or multiple integrals will be taken over
unit tori. 

\begin{thm}
\label{vto}
Let $\f\in(B^2_{\be,1})_+$. Then the function {\em\rf{fun}} has second derivative in the norm and
\bay
\label{fla}
\frac{d^2}{ds^2}\f(T_s)\Big|_{s=t}=2\iiint
(\dg^2\f)(\z,\t,\up)\,d\E_t(\z)\,(R-T)\,d\E_t(\t)\,(R-T)\,d\E_t(\up).
\ey
\end{thm}

\Pf We prove \rf{fla} for $t=0$. For all other $t$ the proof is the same.
By Theorem \ref{per}, 
\begin{align*}
&\frac1t\left(\frac{d}{ds}\big(\f(T_s)\big)\Big|_{s=t}-\frac{d}{ds}\big(\f(T_s)\big)\Big|_{s=0}\right)
\\[.2cm]
=&\frac1t\left(\iint(\dg\f)(\z,\up)\,d\E_t(\z)\,(R-T)\,d\E_t(\up)-
\iint(\dg\f)(\t,\up)\,d\E(\t)\,(R-T)\,d\E(\up)\right)\\[.2cm]
=&\frac1t\left(\iint(\dg\f)(\z,\up)\,d\E_t(\z)\,(R-T)\,d\E_t(\up)-
\iint(\dg\f)(\t,\up)\,d\E(\t)\,(R-T)\,d\E_t(\up)\right)\\[.2cm]
&+\frac1t\left(\iint(\dg\f)(\z,\up)\,d\E(\z)\,(R-T)\,d\E_t(\up)
-\iint(\dg\f)(\z,\t)\,d\E(\z)\,(R-T)\,d\E(\t)\right).
\end{align*}

We have
\begin{align*}
&\iint(\dg\f)(\z,\up)\,d\E_t(\z)\,(R-T)\,d\E_t(\up)-
\iint(\dg\f)(\t,\up)\,d\E(\t)\,(R-T)\,d\E_t(\up)\\[.2cm]
=&\iiint(\dg\f)(\z,\up)\,d\E_t(\z)\,d\E(\t)\,(R-T)\,d\E_t(\up)\\[.2cm]
&-\iiint(\dg\f)(\t,\up)\,d\E_t(\z)\,d\E(\t)\,(R-T)\,d\E_t(\up)\\[.2cm]
=&\iiint(\dg^2\f)(\z,\t,\up)(\z-\t)\,d\E_t(\z)\,d\E(\t)\,(R-T)\,d\E_t(\up)\\[.2cm]
=&\,t\iiint(\dg^2\f)(\z,\t,\up)\,d\E_t(\z)\,(R-T)\,d\E(\t)\,(R-T)\,d\E_t(\up).
\end{align*}

Similarly,
\begin{align*}
&\iint(\dg\f)(\z,\up)\,d\E(\z)\,(R-T)\,d\E_t(\up)
-\iint(\dg\f)(\z,\t)\,d\E(\z)\,(R-T)\,d\E(\t)\\[.2cm]
=&\,t\iiint(\dg^2\f)(\z,\t,\up)\,d\E(\z)(R-T)\,d\E(\t)\,(R-T)\,dE_t(\up).
\end{align*}

Thus
\begin{align*}
&\frac1t\left(\frac{d}{ds}\big(\f(T_s)\big)\Big|_{s=t}-\frac{d}{ds}\big(\f(T_s)\big)\Big|_{s=0}\right)
\\[.2cm]
=&\iiint(\dg^2\f)(\z,\t,\up)\,d\E_t(\z)\,(R-T)\,d\E(\t)\,(R-T)\,d\E_{t}(\up)\\[.2cm]
&+\iiint(\dg^2\f)(\z,\t,\up)\,d\E(\z)\,(R-T)\,d\E(\t)\,(R-T)\,d\E_t(\up).
\end{align*}

As we have mentioned in \S\,2.4, it follows from the results of \cite{Pe8} that  \lb$\dg^2\f\in C_A\hat\otimes C_A\hat\otimes C_A$, i.e., there exist sequences
$\{f_n\}$, $\{g_n\}$, and $\{h_n\}$ in $C_A$ such that
$$
(\dg^2\f)(\z,\t,\up)=\sum_{n\ge0}f_n(\z)g_n(\t)h_n(\up)
$$
and
\bay
\label{kon}
\sum_{n\ge0}\|f_n\|_\be\|g_n\|_\be\|h_n\|_\be<\be.
\ey
Put
$$
Q\df R-T.
$$
We have
$$
\iiint(\dg^2\f)(\z,\t,\up)\,d\E_t(\z)\,Q\,d\E(\t)\,Q\,d\E_{t}(\up)=
\sum_{n\ge0}f_n(T_t)Q\,g_n(T)Q\,h_n(T_t).
$$
It follows from \rf{kon} that
$$
\lim_{t\to0}\sum_{n\ge0}f_n(T_t)Q\,g_n(T)Q\,h_n(T_t)=\sum_{n\ge0}f_n(T)Q\,g_n(T)Q\,h_n(T)
$$
in the operator norm.

Thus
\begin{align*}
\lim_{t\to0}&\iiint(\dg^2\f)(\z,\t,\up)\,d\E_t(\z)\,Q\,d\E(\t)\,Q\,d\E_{t}(\up)\\[.2cm]
=&
\sum_{n\ge0}f_n(T)Q\,g_n(T)Q\,h_n(T)\\[.2cm]
=&\iiint(\dg^2\f)(\z,\t,\up)\,d\E(\z)Q\,d\E(\t)Q\,d\E(\up).
\end{align*}

Similarly,
\begin{align*}
\lim_{t\to0}&\iiint(\dg^2\f)(\z,\t,\up)\,d\E(\z)\,Q\,d\E(\t)\,Q\,d\E_t(\up)\\[.2cm]
=&\iiint(\dg^2\f)(\z,\t,\o)\,d\E(\z)\,Q\,d\E(\t)\,Q\,d\E(\o)
\end{align*}
which proves the result. $\bl$

The same method allows us to prove the following generalization of Theorem \ref{vto}.

\begin{thm}
\label{n}
Suppose that $\f\in\big(B^n_{\be1}\big)_+$. Then the function {\em\rf{fun}} has $n$th derivative
in the norm
$$
\frac{d^n}{ds^n}\f(T_s)\Big|_{s=t}=n!
\underbrace{\int\!\cdots\!\int}_{n+1}(\dg^{n}\f)(\z_1,\cdots,\z_{n+1})\,d\E_t(\z_1)\,(R-T)\cdots(R-T)\,d\E_t(\z_{n+1}).
$$
\end{thm}

\

\section{\bf Differentiability of operator functions in the Hilbert--Schmidt norm}
\setcounter{equation}{0}

\

Suppose that $T$ and $R$ are contractions on Hilbert space such that $T-R\in\bS_2$. We are going to obtain in this section
results on the differentiability of the function \rf{fun} in the Hilbert--Schmidt norm.

\begin{thm}
\label{HSN}
Let $\f$ be a function analytic in $\dd$ such that $\f'\in C_A$. Suppose that $T$ and $R$ are contractions on Hilbert space such that $T-R\in\bS_2$ and and $T_t$ is defined by {\em\rf{Tt}}. Then the function {\em\rf{fun}} is differentiable in the Hilbert--Schmidt norm and
$$
\frac{d}{ds}\f(T_s)\Big|_{s=0}=\iint\frac{\f(\z)-\f(\t)}{\z-\t}\,d\E(\z)\,(R-T)\,d\E(\t),
$$
where $\E$ is a semi-spectral measure of $T$.
\end{thm}

\Pf Let $\f_n$ be the $n$th Ces\'aro mean of the Taylor series of $\f$. Then $\|\f'_n-\f'\|_\be\to0$ as $n\to\be$.
Let $\E_t$ be a semi-spectral measure of $T_t$. We have
\begin{align*}
\frac1s\big(\f(T_s)-\f(T)\big)&=\iint\frac{\f(\z)-\f(\t)}{\z-\t}\,d\E_s(\z)\,(R-T)\,d\E(\t)\\[.2cm]
&=\iint\frac{\f_n(\z)-\f_n(\t)}{\z-\t}\,d\E_s(\z)\,(R-T)\,d\E(\t)\\[.2cm]
+&\iint\frac{(\f-\f_n)(\z)-(\f-\f_n)(\t)}{\z-\t}\,d\E_s(\z)\,(R-T)\,d\E(\t).
\end{align*}

Let $\e>0$. There exists a natural number $N$ such that
\bay
\label{<e}
\sup_{\z,\t}\left|\frac{(\f-\f_n)(\z)-(\f-\f_n)(\t)}{\z-\t}\right|<\e
\ey
whenever $n\ge N$.

Let $n\ge N$. Since $\f_n$ is a polynomial, it is easy to see that
$$
\iint\frac{\f_n(\z)-\f_n(\t)}{\z-\t}\,d\E_s(\z)\,(R-T)\,d\E(\t)
$$
tends to
$$
\iint\frac{\f_n(\z)-\f_n(\t)}{\z-\t}\,d\E(\z)\,(R-T)\,d\E(\t)
$$
in the Hilbert--Schmidt norm. Let $\d>0$ be a number such that
$$
\left|\iint\frac{\f_n(\z)-\f_n(\t)}{\z-\t}\,d\E_s(\z)\,(R-T)\,d\E(\t)-\iint\frac{\f_n(\z)-\f_n(\t)}{\z-\t}\,d\E(\z)\,(R-T)\,d\E(\t)\right|<\e
$$
whenever $s\le\d$.

To conclude the proof, we observe that inequalities \rf{<e} and \rf{fn2} imply
that
$$
\left\|\iint\frac{(\f-\f_n)(\z)-(\f-\f_n)(\t)}{\z-\t}\,d\E_s(\z)\,(R-T)\,d\E(\t)\right\|_{\bS_2}<\e
$$
and
$$
\left\|\iint\frac{(\f-\f_n)(\z)-(\f-\f_n)(\t)}{\z-\t}\,d\E(\z)\,(R-T)\,d\E(\t)\right\|_{\bS_2}<\e
$$
whenever $n\ge N$, and so
$$
\left\|
\frac1s\big(\f(T_s)-\f(T)\big)-\iint\frac{\f(\z)-\f(\t)}{\z-\t}\,d\E(\z)\,(R-T)\,d\E(\t)
\right\|_{\bS_2}<3\e
$$
whenever $s\le\d$. $\bl$

\

\

\noindent
\begin{tabular}{p{8cm}p{14cm}}
Department of Mathematics \\
Michigan State University  \\
East Lansing, Michigan 48824\\
USA
\end{tabular}

\end{document}